\documentclass[12pt]{article}
\usepackage{latexsym}
\usepackage{graphicx}
\usepackage{verbatim}
\usepackage{amssymb}
\usepackage{amsmath}
\title{Higher moments for random multiplicative measures}
\author{K.J. Falconer\\
\small{{\it Mathematical Institute,  
University of St~Andrews, North Haugh, St~Andrews,}} \\
\small{{\it Fife, KY16~9SS, Scotland }}} 
\date{}

\def\bbbn{{\mathbb N}}

\newtheorem{theo}{Theorem}
\newtheorem{prop}[theo]{Proposition}

\newtheorem{lem}[theo]{Lemma}
\newtheorem{cor}[theo]{Corollary}

\newcommand{\ssi}{_{\bf i}} 
\newcommand{\bi}{{\bf i}} 
\newcommand{\bj}{{\bf j}} 
\newcommand{\bu}{{\bf u}} 
\newcommand{\bv}{{\bf v}} 
\newcommand{\bw}{{\bf w}} 
\newcommand{\bp}{{\bf p}} 
\newcommand{\bq}{{\bf q}} 
\newcommand{\ep}{\epsilon} 
\newcommand{\E}{{\sf E}} 
\renewcommand{\P}{{\sf P}} 
\newcommand{\be}{\begin{equation}} 
\newcommand{\ee}{\end{equation}} 
\allowdisplaybreaks

\begin{document}
\maketitle

\begin{abstract} 
We obtain a condition for the $L^q$-convergence of martingales generated by random multiplicative cascade measures for $q>1$ without any self-similarity requirements on the cascades.

\end{abstract}

\medskip

 \section{Introduction}
\setcounter{equation}{0}
\setcounter{theo}{0}

Random multiplicative cascades were introduced as a model for turbulence by Mandelbrot \cite{Man, Man1} in 1974 since when many variants have  been studied. These cascades take the form of a sequence of  random measures $\mu_k$ obtained as the product of random weights indexed by the vertices of an $m$-ary tree, or equivalently as measures on the hierarchy of $m$-ary subintervals of the unit interval. The cascades have the property that the sequence  $\mu_k(A)$ is a martingale for each set $A$, so that $\mu_k$ converges to a random measure $\mu_\infty$. Natural questions about $\mu_\infty$ relate to its non-degeneracy, its moments and its Hausdorff dimension. Such questions were addressed for basic self-similar cascades by Mandelbrot \cite{Man, Man1}, Kahane \cite{Kah} and Kahane and Peyri\`{e}re \cite{Kah,KP}, and  many other aspects and variants of cascades have been developed since, see  \cite{Liu, BM, P} and the references therein. Much of this work concerns cascades where the weights are self-similar, giving rise to a stochastic functional equation which may be solved to establish properties of the limiting measure. Here we obtain conditions for the $q$th moments of the measures $\mu_k(A)$ to be bounded, which implies almost sure and $L^q$-convergence for $q>1$, without the requirement  of self-similarity so that the functional equation methods are not applicable. This was established in the case $1<q \leq 2 $ by Barral and Mandelbrot \cite{BM}, indeed in a more general setting, but the approach here extends to all $q>1$. The analysis required seems significantly more awkward when $q>2$ and in particular when $q$ is not an integer.

To state the main result we need some standard notation.
Random multiplicative cascades are indexed by   a {\it symbolic space} of words formed from the symbols $\{1,2,\ldots,m\}$. 
For $k=0,1,2, \ldots$ let $I_{k}$ be the 
set of all $k$-term words $I_{k} = \{ (i_{1}, i_{2}, \ldots , i_{k}): \, 1 \leq i_{j} \leq 
m \}$, taking $I_{0}$ to consist of the empty word $\emptyset$. 
 We often abbreviate a
word in $I_{k}$ by
${\bf i} = (i_{1}, i_{2}, \ldots , i_{k} ) $ and write $|{\bf i}|=k$ for 
its length.
We let
$I= \cup^{\infty}_{k=0} I_{k} $
denote the set of all finite words with $I_{\infty}$ for the corresponding
set of infinite words, so
$I_{\infty} = \{(i_{1}, i_{2}, \ldots ): 1 \leq i_{j} \leq m \}$. Juxtaposition
of ${\bf i}$ and ${\bf j}$ is written ${\bf ij}$. We write ${\bf 
i}|{k} = (i_{1}, \ldots , i_{k})$ for the {\it curtailment} after $k$ terms
of ${\bf i} = (i_{1}, i_{2}, \ldots ) \in I_{\infty}$, 
or of   ${\bf i} = (i_{1}, \ldots,  i_{k'}) \in I$ if $ k 
\leq k^{\prime}$, with ${\bf i} \preceq {\bf j}$ meaning that ${\bf i}$ is a 
curtailment of ${\bf j}$. If ${\bf i,j} \in I \cup I_{\infty} $  then 
${\bf i}\wedge {\bf j}$ is the maximal word such that both 
${\bf i}\wedge {\bf j} \preceq {\bf i}$ and ${\bf i}\wedge {\bf j} \preceq {\bf j}$. We may topologise $I_{\infty}$ in a natural way by the metric 
$d(\bi,\bj) = m^{-|\bi \wedge \bj |}$ for distinct $\bi,\bj \in 
I_{\infty}$ which makes $I_{\infty}$  into a compact metric space, with
the {\it cylinders}  
$C \ssi = \{\bj \in I_{\infty} : \bi \preceq \bj \}$ for  $\bi \in I$ forming a base 
of open and closed neighbourhoods.
Note that these cylinders are often identified with the hierarchy of $m$-ary subintervals of $[0,1]$ in the natural way, with each  
$C \ssi$ corresponding to a subinterval of length $m^{-|\bi|}$; thus the functions and measures introduced below may be though of as defined on $[0,1]$.

Let $(\Omega, {\cal F}, \P)$ be a probability space. Let  
$\{ W_\bi : \bi \in I\}$ be independent strictly positive random variables, not necessarily identically distributed, with   
$\E(W_\bi) = 1$ for each $\bi \in I$, and set 
$$X_{i_1, \ldots,i_k} = W_{i_1}W_{i_1, i_2}\ldots W_{i_1, \ldots,i_k}.$$
 Let ${\cal F}_k = \sigma( W_\bi : |\bi| \leq k)$ be the $\sigma$-field underlying the random variables indexed by words of length at most $k$. Note that 
$(X_{\bi|k}, {\cal F}_k)$ is a martingale for each $\bi \in I_\infty$; such a family of martingales is termed a $T${\it -martingale}, see  \cite{Kah}.

Let $\mu$ be a given Borel probability measure on $I_\infty$. We may define a sequence of random measures $\mu_k$ on $I_\infty$ by 
\be
\mu_k(A) = \int_A  X_{\bi|k} d\mu(\bi),\label{int}
\ee
for $\mu$-measurable $A$, so in particular for cylinders 
\be
\mu_k(C_\bv) = X_{\bv|k} \mu(C_\bv) \quad \mbox{if } k \leq |\bv|.\nonumber
\ee
It follows that  $(\mu_k(A), {\cal F}_k)$ is a martingale for each Borel set $A$.
By the martingale convergence theorem $\mu_k(A)$ converges  with probability one for each $A$  and $\mu_k$ converges weakly to a measure $\mu_\infty$ on $I_\infty$. 

In particular the sequence of random variables
\be
\mu_k(I_\infty) = \int_{I_\infty}  X_{\bi|k} d\mu(\bi) 
= \sum_{|\bi| = k} X_\bi \mu(C_\bi).\label{mesintv}
\ee
converges almost surely to $\mu_\infty(I_\infty)$.
If $\mu_k(I_\infty)$ is $L^q$-bounded for some $q > 1$, i.e. if
\be
\sup_{k} 
 \E\Big(\Big( \sum_{|\bi|=k}X_\bi \mu(C_\bi) \Big)^q\Big) <\infty,\label{qmom}
 \ee
 then, using Minkowski's inequality, 
\begin{equation}
\limsup_{k \to \infty} \Big(\sum_{|\bi|= k} \E\big((X_\bi \mu(C_\bi))^q\big)\Big)^{1/k}\leq 1.\label{condition2}
\end{equation}
We seek a converse implication, that is a condition such as (\ref{condition2}) for $L^q$ convergence of $\mu_k(I_\infty)$. In the  self-similar case, that is where $(W_{\bi 1},\ldots,W_{\bi m}) $ has the distribution of $(W_{1},\ldots,W_{m}) $ for all $\bi$, (\ref{qmom}) holds if and only if $\E(\sum_{i=1}^m W_i^q) <1$, see, for example, \cite{Kah,KP,Liu}. Here we dispense with any self-similarity assumption.
\begin{theo}\label{mainintA}
Let $X_\bi$ be a  multiplicative random cascade.  Let $q > 1$. Suppose 
\begin{equation}
\limsup_{k \to \infty} \Big(\sum_{|\bi|= k} \E\big((X_\bi \mu(C_\bi))^q\big)\Big)^{1/k} < 1.\label{condition1A}
\end{equation}
Then
\begin{equation}
\limsup_{k \to \infty} 
 \E\Big(\Big( \int X_{\bi |k} d\mu(\bi) \Big)^q\Big) 
= \limsup_{k \to \infty} 
 \E\Big(\Big( \sum_{|\bi|=k}X_\bi \mu(C_\bi) \Big)^q\Big) <\infty . \label{conclusion1}
\end{equation}
 Moreover
 $ \int X_{\bi |k}d\mu(\bi)$ converges almost surely and in   $L^q$ if $q>1$. 
 \end{theo} 
 
Of course, the final conclusion follows immediately from  (\ref{conclusion1}) which says that 
the martingale $ \int X_{\bi |k}d\mu(\bi)$ is $L^q$-bounded. 
 
We will use arguments based on trees and their automorphisms to study the higher moments of multiplicative processes. Tree structures are used in a rather different way to find the $L^q $-dimensions of self-affine sets \cite{Fa5} and the  $L^q$-dimensions of the images of measures under certain Gaussian processes \cite{FX}.

\section{Trees and automorphisms}
\setcounter{equation}{0}
\setcounter{theo}{0}

This section sets out the notation needed relating to the underlying tree structure of the multiplicative cascades.

The integers $m\geq 2$ and $k \geq1$ are fixed throughout this and the next section. For $q>1$ we write $q = n+\epsilon$ where $n$ is an integer and $0 \leq \epsilon <1$. 

We identify the finite words of the symbolic space $I$ with the vertices of the $m$-ary
 rooted tree $T$ with root $\emptyset$ in the natural way. 
We write $T_k$ for the finite rooted tree with vertices $\cup_{l=0}^{k} I_k$. The edges of these trees  join each vertex $\bi \in \cup_{l=0}^{k-1} I_k$  to its $m$  `children'  $\bi 1,\ldots,\bi m$.  Thus the words in $I_k$ are `bottom' vertices of $T_k$.  We will write $T_k$ to denote these trees regarded both as graphs and as sets of vertices, the context making the usage clear.  The estimates in the next section involve automorphisms of the tree $T_k$, regarded as a graph, which induce permutations of the vertices at each level of the tree.

For each vertex $\bv \in T_k$  and $n \geq1$ we write
\be 
S_\bv(n) = \{(\bi_1,\ldots,\bi_n)  \succeq \bv\}\subseteq   ( I_k)^n \label{sdef}
\ee
for the set of all ordered $n$-tuples of $I_k$ (with repetitions allowed)  that are descendents of  $\bv$. 
Let $\mbox{Aut}_\bv$ be the group of automorphisms of the rooted tree $T_k$ that fix $\bv$. Define an equivalence relation $\sim$ on $S_\bv(n)$
by 
\be 
(\bi_1,\ldots,\bi_n)  \sim (\bi_1',\ldots,\bi_n') \mbox{ if there exists } g \in   \mbox{Aut}_\bv  \mbox{ such that  }   g(\bi_r) =  \bi_r'  \mbox{ for  all }  1 \leq r \leq n; \label{equivdef}
\ee
  thus the  equivalence classes  are the orbits of $(I_k)^n$ under $ \mbox{Aut}_\bv $, and we write  $S_\bv(n)/\sim$ for the set of equivalence classes. We write $[J]_\bv$ for the equivalence class containing $J = (\bi_1,\ldots,\bi_n)$.

To work with  the case when $q$ is non-integral  we need to introduce an extra identified point into these orbits.
Let  $\bv \in T_k$ and  let $\bi_1,\ldots,\bi_n $ be (not necessarily distinct) points of $I_k$ such that $\bv \preceq \bi_r$ for all
$r =1,\ldots,n$. Write $T_\bv(\bi_1,\ldots,\bi_n)$  for the minimal subtree of $T_k$ rooted at $\bv$ and containing the points $\{\bi_1,\ldots,\bi_n\}$. 
For each such $\bv \in T_k$  define the ordered set of $(n+1)$-tuples
\be 
S_\bv^+(n)= \{(\bi_1,\ldots,\bi_n;\bp) : (\bi_1,\ldots,\bi_n)  \succeq \bv  \mbox{ and }   \bp \in T_\bv(\bi_1,\ldots,\bi_n) \} \label{sodef}
\ee
(note that in (\ref{sodef})  $\bp$ can be a point  of $T_\bv(\bi_1,\ldots,\bi_n)$ at any level).
We may  define an equivalence relation $\approx$ on each   $S_\bv^+(n)$ by 
\begin{align} 
(\bi_1,\ldots,\bi_n;\bp)  \approx (\bi_1',&\ldots,\bi_n';\bp') \mbox{ if there exists } g \in   \mbox{Aut}_\bv \nonumber \\
& \mbox{ such that  }   g(\bi_r) =  \bi_r'  \mbox{  for all  } 1 \leq r \leq n  \mbox{ and } g(\bp) = \bp'. \label{pequivdef}
\end{align}
We write  $S_\bv^+ (n)/\approx$ for the set of equivalence classes. Observe that, with $ \bp \in T_\bv(\bi_1,\ldots,\bi_n)$, the action of $ \mbox{Aut}_\bv$ on $\bp$ is completely determined by the action of $ \mbox{Aut}_\bv$ on $(\bi_1,\ldots,\bi_n)$.

For notational simplicity, we may omit the subscript when $\bv = \emptyset$, so that
$S(n) \equiv S_\emptyset(n)$, $S^+(n) \equiv S_\emptyset^+(n)$, $[J]= [J]_\emptyset$ and $T(\bi_1,\ldots,\bi_n) \equiv T_\emptyset(\bi_1,\ldots,\bi_n)$.

We require some terminology relating to join sets of elements of  $I_k$. Let $J =  (\bi_1,\ldots,\bi_n)\in (I_k)^n$. The {\it join set} of $J$, denoted by $\bigwedge(J) = \bigwedge(\bi_1,\ldots,\bi_n)$,  is the set of vertices $\{\bv_1,\ldots,\bv_{n-1}\} \in T_k$ consisting of the {\it join points}  $\bi_i\wedge \bi_j $ for all $\bi_i,\bi_j \in I_k$, with $\bw \in \bigwedge(J)$ occurring with {\it multiplicity} $m$ if there are $(m+1)$ distinct $\bi_{i_1},\ldots,\bi_{i_{m+1}}  \in J$ such that  
$\bi_{i_r}\wedge\bi_{i_s} = \bw$ for all $r \neq s$. Note that if two or more of the $\bi_i$ are equal then this common point is automatically a join point with the appropriate multiplicity. Thus $\bigwedge(\bi_1,\ldots,\bi_n)$ consists of the vertices of the subtree $T_\bv(\bi_1,\ldots,\bi_n)$ that have at least two offspring together with any repeated $\bi_i$. Note that the join set of $n$ points always consists of $n-1$ points counting by multiplicity. We write $\wedge^T(\bi_1,\ldots,\bi_n)$ for the `top join point' of $\{\bi_1,\ldots,\bi_n\}$, that is the vertex $\bv$ such that $\bv \preceq \bi_j$ for all $1 \leq j\leq n$ for which $|\bv|$ is greatest.

If $\bj \succeq \bv$ we write $\bj \wedge T_\bv(\bi_1,\ldots,\bi_n)$ for the vertex $\bj| i$ for the largest $i$ such that $\bj| i\in T_\bv(\bi_1,\ldots,\bi_n)$.

The {\it level} of a vertex $\bv \in I$ is just $|\bv|$. Thus the  {\it set of join levels} $L(J)$ of $J \in S_\bv(n)$ is  $\{|\bv_1|,\ldots,|\bv_{n-1}|: \bv_i \in \bigwedge(J)\}$ with levels repeated according to multiplicity. Notice that if $J \sim J'$ then $L(J) = L(J') \equiv L([J]_\bv)$, i.e. the set of levels is constant across each orbit $[J]_\bv$ of $(I_k)^n$.

Similarly, given   $J =(\bi_1,\ldots,\bi_n;\bp) \in S_\bv^+(n)$ we let $L(J)$ be the set of $n-1$ join points of $(\bi_1,\ldots,\bi_n)$ and write $l_0(J) = |\bp|$ for the level of the special point $\bp$. Again this is independent of the choice of $J $ in any orbit in  $S_\bv^+ (n)/\approx$.

\section{The main estimates }
\setcounter{equation}{0}
\setcounter{theo}{0}

This section contains the substance of the proof of Theorem \ref{mainintA} which involves an inductive argument. The induction is significantly more complicated when $q = n +\ep$ is non-integral, that is when $0<\epsilon <1$, when one of the vertices of the underlying tree has to be specifically identified with the `$\ep$' term and we need to work with $S^+(n)$ rather than $S(n)$.

The following identity, which follows from the martingale property and additivity of the measure, will be used repeatedly: for all $\bw \in T_k$ and $|\bw|\leq l \leq k$,
\be
\sum_{|\bv|= l, \bv \succeq \bw} \E\big(X_\bv \mu(C_\bv) \big|{\cal F}_\bw\big)
= \sum_{|\bv|= l, \bv \succeq \bw}X_\bw \mu(C_\bv)
=X_\bw \mu(C_\bw). \label{ident}
\ee
For convenience we define the random variables 
\be
Y_\bv \equiv X_\bv \mu(C_\bv) \qquad (\bv \in I) \label{defY}
\ee
so (\ref{ident}) becomes 
\be
\sum_{|\bv|= l, \bv \succeq \bw} \E\big(Y_\bv \big|{\cal F}_\bw\big)
= Y_\bw \qquad (l \geq |\bw|). \label{ident1}
\ee

Thus  to prove Theorem \ref{mainintA} we must show that
 for $q \geq 1$, 
\begin{equation}
\limsup_{k \to \infty} \Big(\sum_{|\bj|= k} \E\big(Y_\bj^q\big)\Big)^{1/k} < 1 \quad
\mbox{ implies } \quad
\limsup_{k \to \infty} 
 \E\Big(\big( \sum_{|\bi|=k}Y_\bi \big)^q\Big) <\infty.
\label{mainres}
\end{equation}
Note that  (\ref{ident1}) implies that for $0 \leq \ep \leq 1$
\be
\E\Big(\big(\sum_{|\bv|= l, \bv \succeq \bw}Y_\bv\big)^\ep\big|{\cal F}_\bw\Big) 
\leq \Big(\E\sum_{|\bv|= l, \bv \succeq \bw}Y_\bv\big|{\cal F}_\bw\Big) ^\ep = Y_\bw^\ep  \qquad
 (l \geq |\bw|).  \label{epest}
\ee

The strategy of the proof, in the simpler case when $q$ is an integer, is first to estimate the sum of the terms $ \E(Y_{\bi_1}\ldots Y_{\bi_n})$ over the $({\bi_1}\ldots {\bi_n})$ within each equivalence class of $S(n)/\sim$ and then sum these estimates over all the equivalence classes. In the case when $q$ is non-integral, that is when $\epsilon>0$, there is  a further initial stage involving summing within the equivalence classes of $S^+(n)/\approx$.

Note that in   (\ref{intest1}) and below, the product is over the set of levels in a join class. The symbol $[n-1]$ above the product sign merely indicates that there are  $n-1$ terms in this product;  this convention is helpful when keeping track of terms through the proofs.

\begin{prop}\label{integerest}
Let $q =  n+\ep>1$  with $n$ an integer and $0\leq \ep<1$.   Let $J\in S^+(n) $. Then
\be
\sum_{(\bi_1,\ldots,\bi_n;\bp) \in [J]}   \E\Big(Y_{\bi_1}\ldots Y_{\bi_n}
\big(\sum_{\genfrac{}{}{0pt}{1}{|\bj|=k}{ \bj \wedge T(\bi_1,\ldots,\bi_n)=\bp} }Y_\bj\big)^\ep\Big) \leq  \prod_{l \in L(J)}^{[n-1]}
 \Big(\sum_{|\bu|= l} \E\big(Y_\bu^q\big)\Big)^{1/(q-1)}
 \Big(\sum_{|\bu|= l(J)} \E\big(Y_\bu^q\big)\Big)^{\ep/(q-1)}. \label{intest1}
\ee
(When $\ep=0$ the two terms involving $\ep$ disappear and $\bp$ becomes redundant.)
\end{prop}
Proposition \ref{integerest} will follow immediately from the following two lemmas which establish inductive hypotheses that specialize to (\ref{intest1}). The proof of the first lemma, dealing with the case of  $q$  an integer (i.e. with $\epsilon = 0$), is simpler, whilst the proof of the second lemma, for non-integral $q$, both depends on the first result and requires an extension of the approach.

\begin{lem}\label{lemA}
 For all  integers $n \geq 1$, for all $q\geq n$ and all $\bv \in T_k$, if
$J \in S_\bv(n)$  then
\be
  \sum_{(\bi_1,\ldots,\bi_n) \in [J]_\bv} \E\Big(Y_{\bi_1}\ldots Y_{\bi_n}
\big|{\cal F}_\bv\Big)  \leq Y_\bv^{(q-n)/(q-1)} \prod_{l \in L(J)}^{[n-1]}
 \Big(\sum_{|\bu|= l, \bu \succeq \bv} \E\big(Y_\bu^q\big|{\cal F}_\bv\big)\Big)^{1/(q-1)}; \label{intesta}
\ee
\end{lem}

\noindent{\it Proof.} 
\medskip
We obtain (\ref{intesta}) by induction on $n$.

\noindent{\it Start of induction}
If $n=1$,  $\bv \in  T_k$ and $J = (\bj_1) \in S_\bv(1)$ identity (\ref{ident1}) gives
$$
\sum_{ (\bi_1) \in [J]_\bv}   \E\big(Y_{\bi_1}\big|{\cal F}_\bv \big) 
= Y_\bv = Y_\bv^{(q-1)/(q-1)}
$$
which is (\ref{intesta}) when $n = 1$.

\medskip
\noindent{\it The inductive step}
Assume that  for some integer $n_0 \geq 1$ inequality (\ref{intesta})  holds for all  $1 \leq n\leq n_0$, for all $\bv \in T_k$ and all
$J \in S_\bv(n)$.  We establish (\ref{intesta}) when $n = n_0+1$.
Let $\bv \in T_k$ and $J = (\bj_1,\ldots,\bj_n) \in S_\bv(n)$. We divide the argument into two cases.

\medskip
\noindent{\it Case $(a)$}
Assume that   $\bv = \wedge^T(\bj_1,\ldots,\bj_n)$; thus $\bv\preceq \bj_i$ for all $i$ and $\bv$ is itself  the join point of at least two points of $\{\bj_1,\ldots,\bj_n\}$.

If $|\bv|= k$, that is $\bv = \bj_1=\cdots = \bj_n$, then  (\ref{intesta}) is trivially satisfied.

Otherwise $J$ decomposes into   $2\leq r <n$ subsets,
$$J_1= (\bj_1^1,\ldots,\bj_{n_1}^1)\in S_\bv(n_1),  \ldots,\, 
J_{r} = (\bj_1^r,\ldots,\bj_{n_r}^r) \in S_\bv(n_{r}),$$
say, without loss of generality, where $1 \leq n_i \leq n-1$ for each $i$, and
\be
n_1 +\cdots +n_r = n, \label{weightsum}
\ee
and such that each tree 
$T_\bv(\bj_1^{i},\ldots,\bj_{n_{i}}^{i})$ has a distinct {\it single} edge abutting $\bv$. Note that the combinatorics of such a decomposition is preserved under every automorphism in  $\mbox{Aut}_\bv$. We write $L(J_i)\geq 0$ for the set of $(n_i - 1)$ join levels of the trees $T_\bv(\bi_1^{i},\ldots,\bi_{n_{i}}^{i})$  (counted by multiplicity) for $i=1,\ldots,r$.

Using independence conditional  on ${\cal F}_\bv$ and applying the inductive assumption (\ref{intesta}) to $J_1, \ldots J_r$, 
\begin{align*}
&\sum_{(\bi_1,\ldots,\bi_n) \in [J]_\bv}   \E\Big(Y_{\bi_1}\ldots Y_{\bi_n}\big|{\cal F}_\bv\Big) 
\\
& \leq \E\Big(\sum_{(\bi_1^1,\ldots,\bi_{n_1}^1) \in [J_1]_\bv}  
Y_{\bi_1^1}\cdots Y_{\bi_{n_1}^1}\big|{\cal F}_\bv\Big) \times
\cdots \times
\E\Big(\sum_{(\bi_1^{r},\ldots,\bi_{n_{r}}^{r})\in[J_{r}]_\bv} 
Y_{\bi_1^{r}}\cdots Y_{\bi_{n_{r}}^{r}}\big|{\cal F}_\bv\Big)
\\
&\leq 
 Y_\bv^{(q-n_1)/(q-1)} \prod_{l \in L(J_1)}^{[n_1-1]}
 \Big(\sum_{|\bu|= l, \bu \succeq \bv} \E\big(Y_\bu^q\big|{\cal F}_{{\bv}}\big)\Big)^{1/(q-1)}\Big)\times \cdots  \\
& \qquad\times Y_\bv^{(q-n_{r})/(q-1)} 
\prod_{l \in L(J_{r})}^{[n_{r}-1]}   \Big(\sum_{|\bu|= l, \bu \succeq \bv} \E\big(Y_\bu^q\big|{\cal F}_{{\bv}}\big)\Big)^{1/(q-1)}
\\
& =Y_\bv^{(q-n_1-\cdots-n_r)/(q-1)} \big(Y_\bv^q  \big)^{(r-1)/(q-1)}
 \times\prod_{l \in L(J_1)\cup \cdots  \cup  L(J_r)}^{[n_1+ \cdots+ n_r-r]} 
\Big(\sum_{|\bu|= l, \bu \succeq \bv} \E\big(Y_\bu^q\big|{\cal F}_{{\bv}}\big)\Big)^{1/(q-1)}  \\
&= Y_\bv^{(q-n)/(q-1)} 
\prod_{l \in L(J)}^{[n-1]} 
\Big(\sum_{|\bu|= l, \bu \succeq \bv} \E\big(Y_\bu^q\big|{\cal F}_{{\bv}}\big)\Big)^{1/(q-1)},
 \end{align*}
 where we have used (\ref{weightsum}), and incorporated the terms $Y_\bv^q$, taken as a trivial sum over the single vertex $\bv$, in the main product with multiplicity $(r-1)$, to get  (\ref{intesta}) in this case.
 
\medskip
\noindent{\it Case $(b)$}
Now with $\bv \in T_k$ and $J= (\bj_1,\ldots,\bj_n) \in S_\bv(n)$, suppose that  $\bv \preceq \bw_0 = \wedge^T(\bj_1,\ldots,\bj_n)$ and $\bv \neq \bw_0$.
 For each $\bw\succeq \bv$ with 
$|\bw|= l' $ let $g_\bw \in \mbox{Aut}_\bv$ be some automorphism of $T_k$ fixing $\bv$ such that  $g_\bw(\bw_0) = \bw$.
Summing (\ref{intesta}) over each such $\bw$, applying Case (a) to each $g_\bw(J) \in S_\bw(n)$ noting that $L(g_\bw(J)) = L(J)$ and using H\"{o}lder's inequality for the sums and expectations,
\begin{align*}
& \sum_{(\bi_1,\ldots,\bi_n) \in[J]_\bv}\E\Big(Y_{\bi_1}\ldots Y_{\bi_{n}}\big|{\cal F}_\bv\Big) 
\\
& =\E\bigg( \sum_{|\bw| = l', \bw \succeq \bv}\bigg\{ \sum_{(\bi_1,\ldots,\bi_n) \in [g_\bw(J)]_\bw}  \E\Big(Y_{\bi_1}\ldots Y_{\bi_{n}} 
\big|{\cal F}_\bw\Big)\bigg\}\bigg|{\cal F}_\bv\bigg)
\\
& 
\leq\E\bigg( \sum_{|\bw| = l', \bw \succeq \bv}
\bigg\{Y_\bw^{(q-n)/(q-1)} 
\prod_{l \in L(J)}^{[n-1]} 
\Big(\sum_{|\bu|= l, \bu \succeq \bw} \E\big(Y_\bu^q\big|{\cal F}_{{\bw}}\big)\Big)^{1/(q-1)}\bigg\}
\bigg|{\cal F}_\bv\bigg)
\\
& \leq \Big( \E \Big(\sum_{|\bw| = l', \bw \succeq \bv}   Y_\bw  \big|{\cal F}_\bv\Big)\Big)^{(q-n)/(q-1)} \prod_{l \in L(J)}^{[n-1]}
 \Big( \E \Big( \sum_{|\bw| = l', \bw \succeq \bv}\, \sum_{|\bu| = l, \bu \succeq \bw} \E(Y_\bu^q\big|{\cal F}_\bw )\big|{\cal F}_\bv\Big)\Big)^{1/(q-1)}
\\
& \leq Y_\bw  ^{(q-n)/(q-1)} \prod_{l \in L(J)}^{[n-1]}
\Big( \sum_{|\bu| = l, \bu \succeq \bv} \E(Y_\bu^q\big|{\cal F}_\bv )\Big)^{1/(q-1)}
\end{align*} 
using (\ref{ident1}) and the tower property of conditional expectation, giving  (\ref{intesta})  in this case.
$\Box$

The next lemma extends Lemma \ref{lemA} to non-integral $q$. Again conditional independence and H\"{o}lder's inequality are used frequently, but the addition of an extra $\bj \in I_k$ associated with the `$\epsilon$' term significantly complicates the argument.

\begin{lem}\label{lemB}
For all  integers $n \geq 1$, for all $0< \ep<1$, all $q\geq n+\ep$ and all $\bv \in T_k$, if
$J \in S_\bv^+(n)$  then
\begin{align}
 &  \sum_{(\bi_1,\ldots,\bi_n;\bp) \in [J]_\bv} \E\Big(Y_{\bi_1}\ldots Y_{\bi_n}
\big(\sum_{|\bj|=k,\, \bj \wedge T_\bv(\bi_1,\ldots,\bi_n)=\bp}Y_\bj\big)^\ep\big|{\cal F}_\bv\Big) \nonumber\\
&\quad \leq Y_\bv^{(q-n-\ep)/(q-1)} \prod_{l \in L(J)}^{[n-1]}
 \Big(\sum_{|\bu|= l, \bu \succeq \bv} \E\big(Y_\bu^q\big|{\cal F}_\bv\big)\Big)^{1/(q-1)}
 \Big(\sum_{|\bu|= l_{0}(J), \bu \succeq \bv} \E\big(Y_\bu^q\big|{\cal F}_\bv\big)\Big)^{\ep/(q-1)}. \label{intestb}
\end{align}
\end{lem}

\noindent{\it Proof.} 
\medskip
We obtain (\ref{intestb}) by induction on $n$.

\noindent{\it Start of induction}
Let $\bv \in  T_k$ and let $J = (\bj_1;\bq) \in S_\bv^+(1)$  (so  the vertex $\bq$ is on the path from $\bv$ to $\bj_1$). If  $\bj_1= \bq$  it is simple to check  (\ref{intestb}) as the interior sums are over a single term. Otherwise, 
u¤sing  conditional independence,  (\ref{ident1}) and  (\ref{epest}) and then   
H\"{o}lder's inequality,
\begin{align*}
\sum_{ (\bi_1; \bp) \in[J]_\bv}  & \E\Big(Y_{\bi_1}
\big(\sum_{|\bj|=k, \,\bj \wedge T_\bv(\bi_1)=\bp}Y_\bj\big)^\ep\big|{\cal F}_\bv\Big)\\
&  =\sum_{ |\bp|=l_0(J),\, \bp \succeq \bv}  \,
 \E\bigg( \sum_{ \bi_1\succeq \bp} \E\Big(Y_{\bi_1}\big(\sum_{|\bj|=k,\, \bj \wedge T_\bv(\bi_1)=\bp}Y_\bj\big)^\ep\big|{\cal F}_\bp \Big)\Big|{\cal F}_\bv\bigg)
\\
&  \leq\sum_{ |\bp|=l_0(J),\, \bp \succeq \bv}  \,
 \E\bigg( \bigg(\sum_{ \bi_1\succeq \bp} \E\big(Y_{\bi_1}\big|{\cal F}_\bp \big)\bigg)\E\Big(\Big(\sum_{|\bj| = k, \,\bj\succeq \bp}Y_\bj\big)^\ep\big|{\cal F}_\bp \Big)\Big|{\cal F}_\bv\bigg)
\\
&  \leq\sum_{ |\bp|=l_0(J),\, \bp \succeq \bv}  \,
 \E\big(  Y_\bp \,Y_\bp^\ep \big|{\cal F}_\bv\big)
\\
&=\E\Big( \sum_{|\bp | = l_0(J), \bp  \succeq \bv}  Y_\bp ^{(q-1-\ep)/(q-1)} \big( Y_\bp ^q)^{\ep/(q-1)}\big|{\cal F}_\bv\Big)\\
&\leq \Big(\E\Big( \sum_{|\bp | = l_0(J), \bp \succeq \bv}  Y_\bp \big|{\cal F}_\bv\Big)\Big)^{(q-1-\ep)/(q-1)}
 \Big(\E\Big(  \sum_{|\bp | = l_0(J), \bp  \succeq \bv}Y_\bp ^q  \big|{\cal F}_\bv\Big)  \Big)^{\ep/(q-1)}\\
 &= Y_\bv^{(q-1-\ep)/(q-1)}
 \Big(\sum_{|\bp | = l_0(J), \, \bp  \succeq \bv}\E\big( Y_\bp ^q  \big|{\cal F}_\bv\big)  \Big)^{\ep/(q-1)}, 
 \end{align*}
which is  (\ref{intestb}) when $n=1$.

\medskip
\noindent{\it The inductive step}

\noindent Assume that  for some integer $n_0 \geq 1$, inequality (\ref{intestb}) holds for all  $1 \leq n\leq n_0$ and all $J \in S_\bv^+(n)$.  We establish (\ref{intestb}) when $n = n_0+1$.
Let $\bv \in T_k$ and $J = (\bj_1,\ldots,\bj_n;\bq) \in S_\bv^+(n)$. There are three cases.

\medskip
\noindent{\it Case $(a)$}
First assume that   $\bv = \wedge^T(\bj_1,\ldots,\bj_n)$ so $\bv$ is the join point of at least two  of $\{\bj_1,\ldots,\bj_n\}$.

Again, if $|\bv|= k$, that is $\bv = \bj_1=\cdots = \bj_n=\bq$, then  (\ref{intesta}) is straightforward to verify.

Otherwise $J$ decomposes into   $2\leq r <n$ subsets,
$$J_1= (\bj_1^1,\ldots,\bj_{n_1}^1)\in S_\bv(n_1),  \ldots,\, J_{r-1} = (\bj_1^{r-1},\ldots,\bj_{n_{r-1}}^{r-1})\in S_\bv(n_{r-1}) ,\,
J_{r} = (\bj_1^r,\ldots,\bj_{n_r}^r;\bq) \in S^+_\bv(n_{r}),$$
say, without loss of generality, where $1 \leq n_i \leq n-1$ for each $r$, and
\be
n_1 +\cdots +n_r = n, \label{weightsumb}
\ee
and such that each tree 
$T_\bv(\bj_1^{i},\ldots,\bj_{n_{i}}^{i})\, (i=1,\ldots,r)$ has a distinct {\it single} edge abutting $\bv$. We write $L(J_i)$ for the set of $(n_i - 1)$ join levels of the trees $T_\bv(\bi_1^{i},\ldots,\bi_{n_{i}}^{i})$  (counted by multiplicity) for $i=1,\ldots,r$, and $l_0(J_r)  = l_0(J)=
|\bq|$.

Using conditional independence and applying (\ref{intesta}) from Lemma \ref{lemA} and the inductive assumption  (\ref{intestb}) to $J_1, \ldots J_r$, 
\begin{align*}
&\sum_{(\bi_1,\ldots,\bi_n;\bp) \in [J]_\bv}   \E\Big(Y_{\bi_1}\ldots Y_{\bi_n}\big(\sum_{|\bj|=k, \,\bj \wedge T_\bv(\bi_1,\ldots,\bi_n)=\bp}Y_\bj\big)^\ep\big|{\cal F}_\bv\Big) \\
& 
\leq \E\Big(\sum_{(\bi_1^1,\ldots,\bi_{n_1}^1) \in [J_1]_\bv}  
Y_{\bi_1^1}\cdots Y_{\bi_{n_1}^1}\big|{\cal F}_\bv\Big) \times
\cdots \times
\E\Big(\sum_{(\bi_1^{r-1},\ldots,\bi_{n_{r-1}}^{r-1})\in[J_{r-1}]_\bv} 
Y_{\bi_1^{r-1}}\cdots Y_{\bi_{n_{r-1}}^{r-1}}\big|{\cal F}_\bv\Big)
\\
& \hspace{5cm} \times \E\Big(\sum_{(\bi_1^r,\ldots,\bi_{n_r}^r;\bp) \in [J_{r}]_\bv} 
Y_{\bi_1^r}\cdots Y_{\bi_{n_r}^r}\big(\sum_{|\bj|=k,\,\bj \wedge T_\bv(\bi_1^r,\ldots,\bi_{n_r}^r)=\bp}Y_\bj\big)^\ep\big|{\cal F}_\bv\Big)
\\
&\leq 
 Y_\bv^{(q-n_1)/(q-1)} \prod_{l \in L(J_1)}^{[n_1-1]}
 \Big(\sum_{|\bu|= l, \bu \succeq \bv} \E\big(Y_\bu^q\big|{\cal F}_{{\bv}}\big)\Big)^{1/(q-1)}\Big)\times \cdots  \\
  & \qquad \times  Y_\bv^{(q-n_{r-1})/(q-1)} 
\prod_{l \in L(J_{r-1})}^{[n_{r-1}-1]}   \Big(\sum_{|\bu|= l, \bu \succeq \bv} \E\big(Y_\bu^q\big|{\cal F}_{{\bv}}\big)\Big)^{1/(q-1)}\\
&\qquad \times
 Y_\bv^{(q-n_r-\ep)/(q-1)}  
 \prod_{l \in L(J_r)}^{[n_{r}-1]}   \Big(\sum_{|\bu|= l, \bu \succeq \bv} \E\big(Y_\bu^q\big|{\cal F}_{{\bv}}\big)\Big)^{1/(q-1)} \Big(\sum_{|\bu|= l_0(J_r), \bu \succeq \bv} \E\big(Y_\bu^q\big|{\cal F}_{{\bv}}\big)\Big)^{\ep/(q-1)}\\
& =Y_\bv^{(q-n_1-\cdots-n_r -\ep)/(q-1)} \big(Y_\bv^q  \big)^{(r-1)/(q-1)}\\
&\qquad \times\prod_{l \in L(J_1)\cup \cdots  \cup    L(J_r)}^{[n_1+ \cdots+ n_r-r]} 
\Big(\sum_{|\bu|= l, \bu \succeq \bv} \E\big(Y_\bu^q\big|{\cal F}_{{\bv}}\big)\Big)^{1/(q-1)}  \Big(\sum_{|\bu|= l_0(J_r), \bu \succeq \bv}\E\big(Y_\bu^q\big|{\cal F}_{{\bv}}\big)\Big)^{\ep/(q-1)}
 \\
&= Y_\bv^{(q-n-\ep)/(q-1)} 
\prod_{l \in L(J)}^{[n-1]} 
\Big(\sum_{|\bu|= l, \bu \succeq \bv} \E\big(Y_\bu^q\big|{\cal F}_{{\bv}}\big)\Big)^{1/(q-1)}  \Big(\sum_{|\bu|= l_0(J), \bu \succeq \bv}\E\big(Y_\bu^q\big|{\cal F}_{{\bv}}\big)\Big)^{\ep/(q-1)},
 \end{align*}
using (\ref{weightsumb}), and incorporating the terms $Y_\bv^q$ in the main product with multiplicity $(r-1)$ to get  (\ref{intestb}).

\medskip
\noindent{\it Case $(b)$}
Now with $\bv \in T_k$ and $J= (\bj_1,\ldots,\bj_n;\bq) \in S^+_\bv(n)$, suppose that  $\bv \preceq \bw_0 = \wedge^T(\bj_1,\ldots,\bj_n)$ and $\bv \neq \bw_0$.
Also suppose that $\bq \succeq \bw_0$ so that   
$J \in S^+_{\bw_0}(n)$, and let $l'= |\bw_0|>|\bv| $. For each $\bw\succeq \bv$ with 
$|\bw|= l' $ let $g_\bw \in \mbox{Aut}_\bv$ be some automorphism of $T_k$ fixing $\bv$ such that  $g_\bw(\bw_0) = \bw$.
By Case (a) (\ref{intestb}) is valid with $\bv$ replaced by each such $\bw$ in turn, so summing over $\bw$   and using H\"{o}lder's inequality,
\begin{align*}
& \sum_{(\bi_1,\ldots,\bi_n;\bp) \in[J]_\bv}\E\Big(Y_{\bi_1}\ldots Y_{\bi_{n}} \big(\sum_{|\bj|=k,\,\bj \wedge T_\bv(\bi_1,\ldots,\bi_n)=\bp}Y_\bj\big)^\ep\big|{\cal F}_\bv\Big) \\
& =
\E\bigg( \sum_{|\bw| = l', \bw \succeq \bv}\bigg\{ \sum_{(\bi_1,\ldots,\bi_n;\bp) \in [g_\bw(J)]_\bw}  \E\Big(Y_{\bi_1}\ldots Y_{\bi_{n}} 
\big(\sum_{|\bj|=k,\,\bj \wedge T_\bw(\bi_1,\ldots,\bi_n)=\bp}Y_\bj\big)^\ep\big|{\cal F}_\bw\Big)\bigg\}\bigg|{\cal F}_\bv\bigg)\\ 
& \leq \E \bigg(\sum_{|\bw| = l', \bw \succeq \bv} \bigg\{ Y_\bw^{(q-n-\ep)/(q-1)} \prod_{l \in L(J)}^{[n-1]}
 \Big(\sum_{|\bu|= l, \bu \succeq \bw} \E\big(Y_\bu^q\big|{\cal F}_\bw\big)\Big)^{1/(q-1)}
  \Big(\sum_{|\bu|= l_0, \bu \succeq \bw} \E\big(Y_\bu^q\big|{\cal F}_\bw\big)\Big)^{\ep/(q-1)}\bigg\}\bigg|{\cal F}_\bv\bigg) \\
 & \leq \Big( \E \Big(\sum_{|\bw| = l', \bw \succeq \bv}   Y_\bw  \big|{\cal F}_\bv\Big)\Big)^{(q-n-\ep)/(q-1)} \prod_{l \in L(J)}^{[n-1]}
 \Big( \E \Big( \sum_{|\bw| = l', \bw \succeq \bv}\, \sum_{|\bu| = l, \bu \succeq \bw} \E(Y_\bu^q\big|{\cal F}_\bw )\big|{\cal F}_\bv\Big)\Big)^{1/(q-1)}\\
&\hspace{6cm}
\times \Big( \E \Big( \sum_{|\bw| = l', \bw \succeq \bv}\, \sum_{|\bu| = l_0(J), \bu \succeq \bw} \E(Y_\bu^q\big|{\cal F}_\bw )\big|{\cal F}_\bv\Big)\Big)^{\ep/(q-1)}\\
&  = Y_\bv^{(q-n-\ep)/(q-1)} \prod_{l \in L(J)}^{[n-1]}
 \Big(  \sum_{|\bu| = l, \bu \succeq \bv} \E(Y_\bu^q\big|{\cal F}_\bv )\Big)^{1/(q-1)} 
 \Big(  \sum_{|\bu| = l_0(J), \bu \succeq \bv} \E(Y_\bu^q\big|{\cal F}_\bv\big) \Big)^{\ep/(q-1)},
\end{align*}
using (\ref{ident1}) to get  (\ref{intestb})  in this case.

\medskip
\noindent{\it Case $(c)$}
With $\bv \in T_k$ and $J= (\bj_1,\ldots,\bj_n;\bq) \in S^+_\bv(n)$ write $J^- = (\bj_1,\ldots,\bj_n) \in S_\bv(n)$.
As in Case (b) suppose that  $\bv \preceq \bw_0 = \wedge^T(\bj_1,\ldots,\bj_n)$ and $\bv \neq \bw_0$, but now with  $\bv \preceq \bq \preceq \bw_0$ ($\bq \neq \bw_0$) so that $\bq$ lies on the path joining $\bv$ to $\bw_0$.  
For each $\bp\succeq \bv$ with $|\bp|=|\bq| = l_0(J)$ let $g_\bp \in \mbox{Aut}_\bv$ be some tree automorphism fixing $\bv$ with  $g_\bp(\bq) = \bp$.
Spliting the sum over $ \bp \succeq \bv$, using conditional independence, (\ref{intesta}) and (\ref{epest}), and again applying H\"{o}lder's inequality,
\begin{align*}
  & \sum_{(\bi_1,\ldots,\bi_n;\bp) \in[J]_\bv}\E\Big(Y_{\bi_1}\ldots Y_{\bi_{n}} \big(\sum_{|\bj|=k,\,\bj \wedge T_\bv(\bi_1,\ldots,\bi_n)=\bp}Y_\bj\big)^\ep\big|{\cal F}_\bv\Big) \\
& =
\E\bigg( \sum_{|\bp| = l_0(J), \bp \succeq \bv}\bigg\{ \sum_{(\bi_1,\ldots,\bi_n) \in [g_\bp(J^-)]_\bp}  \E\Big(Y_{\bi_1}\ldots Y_{\bi_{n}} 
\big(\sum_{|\bj|=k,\,\bj \wedge T_\bv(\bi_1,\ldots,\bi_n)=\bp}Y_\bj\big)^\ep\big|{\cal F}_\bp\Big)\bigg\}\bigg|{\cal F}_\bv\bigg)\nonumber \\ 
& \leq
\E\bigg( \sum_{|\bp| = l_0(J), \bp \succeq \bv}\bigg\{\Big( \sum_{(\bi_1,\ldots,\bi_n) \in [g_\bp(J^-)]_\bp}  \E\big(Y_{\bi_1}\ldots Y_{\bi_{n}} \big| {\cal F}_\bp  \big)\Big)
\Big(\E\big(\sum_{|\bj|=k,\,\bj \wedge T_\bv(\bi_1,\ldots,\bi_n)=\bp}Y_\bj\big)^\ep\big|{\cal F}_\bp\Big)\bigg\}\bigg|{\cal F}_\bv\bigg)\nonumber \\ 
& \leq \E\bigg( \sum_{|\bp| = l_0(J), \bp \succeq \bv} \bigg\{ Y_\bp^{(q-n)/(q-1)} \prod_{l \in L(J)}^{[n-1]}
 \Big(\sum_{|\bu|= l, \bu \succeq \bp} \E\big(Y_\bu^q\big|{\cal F}_\bp\big)\Big)^{1/(q-1)}\,\,
 Y_\bp^\ep\bigg\}\bigg|{\cal F}_\bv\bigg) \\
& \leq \E\bigg( \sum_{|\bp| = l_0(J), \bp \succeq \bv} \bigg\{ Y_\bp^{(q-n-\ep)/(q-1)} \prod_{l \in L(J)}^{[n-1]}
 \Big(\sum_{|\bu|= l, \bu \succeq \bp} \E\big(Y_\bu^q\big|{\cal F}_\bp\big)\Big)^{1/(q-1)}\,\,
( Y_\bp^q)^{\ep/(q-1)}\bigg\}\bigg|{\cal F}_\bv\bigg) \\
& \leq \E\Big( \sum_{|\bp| = l_0(J), \bp \succeq \bv}  Y_\bp\Big|{\cal F}_\bv\Big)^{(q-n-\ep)/(q-1)} \prod_{l \in L(J)}^{[n-1]}
\E\Big( \sum_{|\bp| = l_0(J), \bp \succeq \bv}  \Big(\sum_{|\bu|= l, \bu \succeq \bp} \E\big(Y_\bu^q\big|{\cal F}_\bp\big)\Big)\Big|{\cal F}_\bv\Big)^{1/(q-1)}\\
&
\hspace{3cm} \times\E\Big( \sum_{|\bp| = l_0(J), \bp \succeq \bv}  Y_\bp^q\Big|{\cal F}_\bv\Big)^{\ep/(q-1)} \\
 & = Y_\bv^{(q-n-\ep)/(q-1)} \prod_{l \in L(J)}^{[n-1]}                         
 \Big(  \sum_{|\bu| = l, \bu \succeq \bv} \E(Y_\bu^q\big|{\cal F}_\bv )\Big)^{1/(q-1)} 
 \Big(  \sum_{|\bu| = l_0, \bu \succeq \bv} \E(Y_\bu^q\big|{\cal F}_\bv\big) \Big)^{\ep/(q-1)},
  \end{align*}
giving  (\ref{intestb}).

This completes the inductive step and the proof of the lemma. $\Box$

\medskip

Proposition \ref{integerest} follows immediately from Lemmas \ref{lemA} and \ref{lemB}  on taking $\bv = \emptyset$.

\medskip

For the case of non-integral $q$ we now  sum  inequality (\ref{intest1}) over all $\bp \in T(J)$   where now  $J \in S(n)$. Recall that  $T(J) \equiv  T_\emptyset(\bj_1,\ldots,\bj_n)$  for $J = (\bj_1,\ldots,\bj_n) \in S(n)$ .

\begin{cor}\label{corollary1}
Let $J \in S(n)$. If $q=n$ is an integer then
\be
\sum_{(\bi_1,\ldots,\bi_n) \in [J]} \E\Big( Y_{\bi_1}\ldots Y_{\bi_n}\Big)
\leq  \prod_{l \in L(J)}^{[n-1]}
 \Big(\sum_{|\bu|= l} \E\big(Y_\bu^q\big)\Big)^{1/(q-1)}. \label{coreq1}
\ee
If $q =  n+\ep>1$  with $n$ an integer and $0< \ep<1$ then
\be
\sum_{(\bi_1,\ldots,\bi_n) \in [J]} \E\Big( Y_{\bi_1}\ldots Y_{\bi_n}\big(\sum_{|\bj|=k}Y_\bj\big)^\ep\Big)
\leq  \sum_{\bq \in T(J)}\prod_{l \in L(J)}^{[n-1]}
 \Big(\sum_{|\bu|= l} \E\big(Y_\bu^q\big)\Big)^{1/(q-1)}
\Big(\sum_{|\bu|= |\bq|} \E\big(Y_\bu^q\big)\Big)^{\ep/(q-1)}. \label{coreq2}
\ee
\end{cor}

\noindent{\it Proof.} 
Inequality (\ref{coreq1}) is just Proposition \ref{integerest} with $\epsilon = 0$.

For (\ref{coreq2}), if $J \in S(n)$ and $\bq \in T(J)$, write $J_\bq = (\bj_1,\ldots,\bj_n,q) \in S^+(n)$. With $0<\epsilon <1$, breaking up the sum and using that $(\sum a_i)^\ep \leq \sum a_i^\ep$ for  $a_i \geq 0$, 
\begin{align*}
\sum_{(\bi_1,\ldots,\bi_n) \in [J]} \E\Big(  Y_{\bi_1}\ldots Y_{\bi_n}\big(\sum_{|\bj|=k}Y_\bj\big)^\ep\Big)
&=  \sum_{(\bi_1,\ldots,\bi_n)  \in [J] }   \E\Big(Y_{\bi_1}\ldots Y_{\bi_n}
  \Big(\sum_{\bp \in T(\bi_1,\ldots,\bi_n)}\,\sum_{|\bj|=k,\,\bj \wedge T(\bi_1,\ldots,\bi_n)=\bp} Y_\bj\Big)^{\ep}\Big)\nonumber\\
&\leq  \sum_{(\bi_1,\ldots,\bi_n)  \in [J] }   \E\Big(Y_{\bi_1}\ldots Y_{\bi_n}
 \sum_{\bp \in T(\bi_1,\ldots,\bi_n)} \big(\sum_{|\bj|=k,\,\bj \wedge T(\bi_1,\ldots,\bi_n)=\bp} Y_\bj\big)^{\ep}\Big)\nonumber\\
 &=  \sum_{(\bi_1,\ldots,\bi_n)  \in [J]}   \sum_{\bp \in T(\bi_1,\ldots,\bi_n)} \E\Big(Y_{\bi_1}\ldots Y_{\bi_n}
\big(\sum_{|\bj|=k,\, \bj\wedge T(\bi_1,\ldots,\bi_n)=\bp} Y_\bj\big)^{\ep}\Big)\nonumber\\
 &\leq    \sum_{\bq \in T(J) }\,\, \sum_{(\bi_1,\ldots,\bi_n,\bp)  \in  [J_\bq]}  \E\Big(Y_{\bi_1}\ldots Y_{\bi_n}
\big(\sum_{|\bj|=k,\, \bj\wedge T(\bi_1,\ldots,\bi_n)=\bp} Y_\bj\big)^{\ep}\Big)\nonumber\\
&\leq    \sum_{\bq \in T(J)}\,\, \prod_{l \in L(J_\bq)}^{[n-1]}
 \Big(\sum_{|\bu|= l} \E\big(Y_\bu^q\big)\Big)^{1/(q-1)}
 \Big(\sum_{|\bu|= l_0(J_\bq)} \E\big(Y_\bu^q\big)\Big)^{\ep/(q-1)}, 
\end{align*}
using Proposition \ref{integerest}, and this is just (\ref{coreq2}) since $L(J_\bq) = L(J)$ and $l_0(J_q) = |\bq|$.
$\Box$

\section{Completion of the proof}
\setcounter{equation}{0}
\setcounter{theo}{0}

Finally, we  have to sum (\ref{coreq1}) and  (\ref{coreq2}) over all equivalence classes $[J]$ of $ S_\bv(n)$ under $\sim$, and to do this we need  to bound the number of equivalence classes with given sets of levels. Let $0 \leq l_1\leq \cdots \leq l_{n} \leq k$ and $0\leq l \leq k$ be (not necessarily distinct) levels. Write
$$
N(l_1,\ldots, l_{n-1}) =\, \#\big\{[J] \in S(n)/\sim  \mbox{ such that }
L([J]) = \{l_1,\ldots, l_{n-1}\} \big\}
$$
and
$$
N^+(l_1,\ldots, l_{n-1};l) =\, \#\big\{[J] \in S^+(n)/\approx \mbox{ such that }
L([J]) = \{l_1,\ldots, l_{n-1}\}, l([J])= l  \big\}. 
$$

\begin{lem}\label{count}
Let $0<\lambda<1$. For $n \in \bbbn$
\be
\sum_{0 \leq l_1\leq \cdots  \leq l_{n-1} \leq k}N(l_1,\ldots,l_{n-1})  \lambda^{l_1+\cdots +l_{n-1}}  \leq M <\infty,\label{count1}
\ee
and  for $n \in \bbbn$ and $\ep>0$
\be
\sum_{0 \leq l_1\leq \cdots  \leq l_{n-1} \leq k, \,0\leq l \leq k}N^+(l_1,\ldots,l_{n-1};l)  \lambda^{l_1+\cdots +l_{n-1}+\ep l}  \leq M^+ <\infty,\label{count2}
\ee
where the bounds $M$ and $M^+$are independent of $k$.
\end{lem}
{\it Proof.}   
For $n \geq 1$, every $J \in S(n+1)$ with $L(J)=\{ l_1,\ldots, l_{n-1},l_n\}$ where $0 \leq l_1\leq \cdots \leq l_{n}$ may be obtained by adjoining a vertex $\bj \in I_k$ to some $J^- \in S(n)$ where
$L(J^-) = \{ l_1,\ldots, l_{n-1}\}$. For each  $J^- \in S(n)$ such a vertex may be adjoined so that the additional join level is $l_n$ in at most $n$ non-equivalent ways under $\mbox{Aut}_\emptyset$.
Thus
$N(l_1,\ldots, l_{n}) \leq n N (l_1,\ldots, l_{n-1})$ and $N (l_1)=1$, so
$N(l_1,\ldots, l_n) \leq  n!$.
Thus, for $0< \lambda<1$, 
 \begin{align*}
\sum_{0 \leq l_1\leq \cdots  \leq l_{n-1} \leq k}N(l_1,\ldots,l_{n-1})  \lambda^{l_1+\cdots +l_{n-1}}  
&\leq (n-1)!\sum_{0 \leq l_1\leq \cdots \leq l_{n-1}\leq k}  \lambda^{ l_1+\cdots +l_{n-1}}\\
&\leq
(n-1)! \sum_{r=0}^{\infty}P(r) \lambda^{r}  \equiv M 
\end{align*}
where $P(r)$ is the number of distinct ways of partitioning the integer $r$ into a sum of $n$ integers
$r = l_1 + \cdots + l_{n}$ where $0 \leq l_1\leq \cdots \leq l_n$. Since $P(r)$ is polynomially bounded (trivially $P(r) \leq (r+1)^{n-1}$) the series is convergent.

Furthermore, for each $J\in S(n)$ there are at most $n$  ways of choosing a vertex $\bq \in T(J)$ at level $l$, so
$$N^+(l_1,\ldots, l_{n-1};l) \leq n N (l_1,\ldots, l_{n-1})  \leq  n(n-1)! = n!.$$ Thus, for $0< \lambda<1$ and $\ep>0$,
 \begin{align*}
\sum_{0 \leq l_1\leq \cdots \leq l_{n-1}\leq k, \,0\leq l \leq k}N^+(l_1,\ldots,l_{n-1};l)  \lambda^{l_1+\cdots +l_n+\ep l}
&\leq
n!\sum_{0 \leq l_1\leq \cdots \leq l_{n-1}\leq k, \,0\leq l \leq k}  \lambda^{ l_1+\cdots +l_{n-1}+\ep l}\nonumber\\
&\leq
n!\, \sum_{r=1}^{\infty}\frac{P(r) \lambda^{r}}{1-\lambda^{\epsilon}}  \equiv M^+.
\end{align*}
$\Box$

\medskip
To get the final conclusion we now sum the inequalities of  Corollary \ref{corollary1} over all equivalence classes $[J] \in S(n)/\sim$ using Lemma \ref{count}  to bound the sums.

\begin{prop}\label{mainint}
Let $q \geq 1$. Suppose 
\begin{equation}
\limsup_{k \to \infty} \Big(\sum_{|\bi|= k} \E\big(Y_\bi^q\big)\Big)^{1/k} < 1.\label{condition}
\end{equation}
Then
\begin{equation}
\limsup_{k \to \infty} 
 \E\Big(\Big( \sum_{|\bi|=k}Y_\bi \Big)^q\Big) <\infty . \label{conclusion}
\end{equation}
\end{prop} 
{\it Proof.}  
The case of $q=1$ is trivial, so assume that $q>1$. From (\ref{condition}) there are numbers $c>0$ and $0<\lambda <1$ such that
\be
\sum_{|\bj|= k} \E\big(Y_\bj^q\big) \leq c \lambda^{(q-1)k} \quad (k= 1,2,\ldots),
\label{condition1}
\ee
choosing $\lambda$ so that $\lambda^{(q-1)}$ is sufficiently close to $1$. In the case of  $q= n+\epsilon$, where $n$ is an integer and $0< \ep<1$, rearranging the summation and using (\ref{coreq2}) gives

\begin{align*}
 \E\Big(\big( \sum_{|\bi|=k}Y_\bi \big)^q\Big)
 &= 
  \E\Big(\big( \sum_{|\bi_1|=k}Y_{\bi_1} \big)\cdots
  \big( \sum_{|\bi_n |=k}Y_{\bi_n} \big)\big( \sum_{|\bj|=k}Y_\bj \big)^\ep\Big)\\
  &= 
\sum_{(\bi_1,\cdots,\bi_n) \in S(n)} \E\Big(  Y_{\bi_1}\cdots Y_{\bi_n}\big(\sum_{|\bj|=k}Y_\bj\big)^\ep\Big)\\
 &= 
\sum_{[J]\in S(n)/\sim}\, 
\sum_{(\bi_1,\ldots,\bi_n) \in [J]}\E\Big(  Y_{\bi_1}\cdots Y_{\bi_n}\big(\sum_{|\bj|=k}Y_\bj\big)^\ep\Big)\\
 &= \sum_{[J]\in S(n)/\sim}\, 
 \sum_{\bq \in T(J)}\,\, \prod_{l \in L(J)}^{[n-1]}
 \Big(\sum_{|\bu|= l} \E\big(Y_\bu^q\big)\Big)^{1/(q-1)}
 \Big(\sum_{|\bu|= |\bq|} \E\big(Y_\bu^q\big)\Big)^{\ep/(q-1)} \\
 &= \sum_{[J]\in S(n)/\sim}\, 
 \sum_{\bq \in T(J)}\,\, \prod_{l \in L(J)}^{[n-1]}
 \big(c\lambda^l\big)
 \big(c\lambda^{ |\bq|}\big)^{\ep} \\
&=c^{(n-1+\ep)}
 \sum_{(J,\bq)\in S^+(n)/\approx}\,\, \prod_{l \in L(J)}^{[n-1]}
 \big(\lambda^l\big)
 \big(\lambda^{|\bq|}\big)^{\ep} \\
 &= c \sum_{0 \leq l_1\leq \cdots  \leq l_{n-1} \leq k, \,0\leq l \leq k}N^+(l_1,\ldots,l_{n-1};l)  \lambda^{(l_1+\cdots +l_{n-1}+\ep l)}\\
 &  \leq M^+ <\infty,
\end{align*}
by (\ref{count2}).

The case where $q=n$ is an integer is similar but shorter, using (\ref{coreq1}) and (\ref{count1}).
$\Box$.

\medskip

{\it Proof of Theorem 1.1}  
Setting
$$
Y_\bi = X_\bi \mu(C_\bv) \qquad (\bi \in I)
$$
in (\ref{condition}) and  (\ref{conclusion}),  Proposition \ref{mainint} gives (\ref{conclusion1}) directly.
Moreover, (\ref{conclusion1})  says that
$( \int X_{\bi |k} d\mu(\bi)  |{ \cal F}_k )$
is an $L^q$-bounded martingale, so the martingale  convergence theorem implies  almost sure convergence and convergence in $L^q$ for $q>1$.
$\Box$


\medskip



\begin{thebibliography}{abc-12}
 \bibitem{BM} 
 J. Barral and B. Mandelbrot. Introduction to infinite products of random independent functions (Random multiplicative multifractal measures), in ÒFractal Geometry and Applications: A Jubilee of Benoõt MandelbrotÓ. M. L. Lapidus and M. van Frankenhuysen eds. {\em Proc. Symp. in Pure Math.} {\bf 72} AMS, Providence, RI (2004): Part I, pp 3--16; Part II, pp 17--52 ; Part III, pp 53--90.
  \bibitem{Fa5} K.J.~Falconer. Generalized dimensions of measures on
   almost self-affine sets, {\em Nonlinearity }{\bf 23}(2010), 1047--1069.
 \bibitem{FX} K.J. Falconer and Yimin Xiao. 
Generalized dimensions of images of measures under Gaussian processes, 
arXiv:1212.2383 (2012).
 \bibitem{Kah} J.-P. Kahane. Sur le mod\`{e}le de turbulence de Beno\^{i}t Mandelbrot. {\em C.R. Acad.
Sci. Paris} {\bf 278}(1974), 567--569.
 \bibitem{KP} J.-P. Kahane, J. Peyri\`{e}re. Sur certaines martingales de B. Mandelbrot. {\em Adv.
Math.} {\bf 22}(1976), 131--145.
 \bibitem{Liu} Quansheng Liu. On generalized multiplicative cascades. {\em Stoch. Proc. Appl.} {\bf 86}(2000), 263-286.
 \bibitem{P} J. Peyri\`{e}re. Recent results on Mandelbrot multiplicative cascades. {\em Progress in Probab.} {\bf 46}(2000), Birkh\"{a}user, pp. 147--159.
 \bibitem{Man} B.B Mandelbrot. Multiplications al\'{e}atores it\'{e}r\'{e}es et distributions invariantes par moyennes pond\'{e}r\'{e}es. {\em C. R. Acad. Sci. Paris} {\bf  278}(1974), 289--292 and 355--358.
 \bibitem{Man1} B.B Mandelbrot. Intermittent turbulence in self-similar cascades: divergence of higher moments and dimension of the carrier. {\em J. Fluid Mech.} {\bf  62}(1974), 331--358.


\end{thebibliography}
\end{document}